\documentclass[12pt]{article}
\usepackage{epsfig}
\usepackage{color}
\usepackage{epstopdf}
\usepackage[mathscr]{eucal}
\usepackage{amsmath}
\usepackage{amssymb}
\usepackage{graphics}
\usepackage{psfrag}

\newtheorem{thm}{Theorem}[section]

\newtheorem{prop}[thm]{Proposition} 
\newtheorem{rmk}[thm]{Remark} 

\newtheorem{examp}[thm]{Example}
\newtheorem{hyp}[thm]{Hypothesis}

\newcommand{\SEone}{{\bf SE}(1)}

\newcommand{\proof}{\noindent {\bf Proof} \hspace{0.2in}} 
\newcommand{\qed}{\hfill\mbox{\raggedright\rule{.07in}{.1in}}
  \vspace{1ex}} 

\newcommand{\Section}[1]{\section{#1} \setcounter{equation}{0}}

\title{Forced translational symmetry-breaking for abstract evolution equations: the organizing center for blocking of travelling waves}
 
\author{Victor G. LeBlanc\\Department of
  Mathematics and Statistics\\University of Ottawa\\Ottawa, 
ON K1N 6N5\\CANADA
\and
Christian Roy\\
Department of
  Mathematics and Statistics\\University of Ottawa\\Ottawa, 
ON K1N 6N5\\CANADA}
\date{\today}

\begin{document}

\maketitle

\pagebreak

\begin{abstract}
We consider two parameter families of differential equations on a Banach space $X$, where the parameters $c$ and $\varepsilon$ are such that:
\begin{itemize}
\item when $\varepsilon=0$, the differential equations are symmetric under the action of the group of one-dimensional translations $\SEone$ acting on $X$, whereas when $\varepsilon\neq 0$, this translation symmetry is broken,
\item when $\varepsilon=0$, the symmetric differential equations admit a smooth family of relative equilibria (travelling waves) parametrized by the drift speed $c$, with $c=0$ corresponding to steady-states.
\end{itemize}
Under certain hypotheses on the differential equations and on the Banach space $X$, we use the center manifold theorem of Sandstede, Scheel and Wulff \cite{SSW} to study the effects of the symmetry-breaking perturbation on the above family of relative equilibria.  In particular,  we show that the phenomenon commonly referred to as {\em propagation failure}, or {\em wave blocking} occurs in a cone in the $(c,\varepsilon)$ parameter space which emanates from the point $(c,\varepsilon)=(0,0)$.

We also discuss how our methods can be adapted to perturbations of parameter-independent differential equations (such as the Fisher-KPP) which admit families of relative equilibria parametrized by drift speed.
\end{abstract}

\pagebreak
\Section{Introduction}
\renewcommand{\thefootnote}{\ensuremath{\fnsymbol{footnote}}}

Travelling waves are an important class of solutions for many types of systems which are modelled using reaction-diffusion partial differential equations \cite{Ei,FM,GK,IM,Jones,KeenerSneyd,LewisKeener,TYBN}, integro-differential equations \cite{BNPR,KFB,SRLC} or delayed partial differential equations \cite{BMR,FT,GXZY,ReyMackey}.  In biological systems, travelling waves frequently represent electrical pulses, fronts, backs or wave trains which propagate along one-dimensional networks of excitable cells \cite{KeenerSneyd,Murray}.    They have also been observed in chemical reactions, combustion theory \cite{VV} and nonlinear optics \cite{AA}.   In fact, the literature on travelling waves, both theoretical and applications, is quite extensive.  We will however single out the excellent review articles by Sandstede \cite{Sandstede} and by Xin \cite{Xin}.

For nonlinear models of propagation in spatially extended systems, a common modelling hypothesis or simplifying assumption is to assume that the underlying medium of propagation of the wave is homogeneous.  However, this hypothesis is seldom reasonable in many practical applications.  Indeed, biological media typically contain many sources of inhomogeneity (e.g. eschemic tissue, gap junctions), which may have effects on the propagation of the wave through the medium (i.e. acceleration, deceleration, reflection, or blocking propagation) \cite{BMR,IM,LewisKeener,TYBN,Xin}.   It is this interplay between inhomogeneity of the medium and the propagation properties of travelling waves which we seek to elucidate in this paper, from {\em a model-independent} point of view.

Mathematically speaking, travelling waves are typically a consequence of translation-invariance (homogeneity) of the underlying mathematical model.  To illustrate this point, consider the following well-known bistable reaction-diffusion equation (see e.g. \cite{KeenerSneyd})
\begin{equation}
\frac{\partial u}{\partial t}=\frac{\partial^2 u}{\partial x^2}+u(1-u)(u-\alpha)
\label{bis}
\end{equation}
where $u(x,t)$ is a function of $x\in\mathbb{R}$ and $t\geq 0$, and $0<\alpha<1$ is a real parameter.  Equation (\ref{bis}) is invariant under the change of variables
\begin{equation}
v(x,t)=u(x+\gamma,t)
\label{e1actfunc}
\end{equation}
for any $\gamma\in\mathbb{R}$.   Using the language of group-equivariant dynamical systems \cite{GSS2}, we say that (\ref{e1actfunc}) defines an action of the additive group of real numbers $\SEone$ on a suitable space of functions, and that the right-hand side of equation (\ref{bis}) is equivariant under this action.  One consequence of this symmetry property (see e.g. \cite{GSS2}) is that solutions of (\ref{bis}) are mapped to other solutions by elements of the group $\SEone$.  A travelling wave solution to (\ref{bis}) is a solution such that its time orbit is contained in its $\SEone$-orbit.  Such a solution is called a {\em relative equilibrium}.  

To find relative equilibria in (\ref{bis}), we assume the ansatz
\begin{equation}
u(x,t)=u^*(c)(x+ct),
\label{twansatz}
\end{equation}
which reduces (\ref{bis}) into an ordinary differential equation for the function $u^*(c)(\eta)$:
\begin{equation}
c\frac{du^{*}(c)}{d\eta}=\frac{d^2 u^{*}(c)}{d\eta^2}+u^*(c)(1-u^*(c))(u^*(c)-\alpha).
\label{bisode}
\end{equation}
One then uses techniques from phase-plane analysis to seek values of $c$ for which (\ref{bisode}) admits a heteroclinic orbit connecting the equilibrium points 
${\displaystyle \left(u^*(c),\frac{d u^*(c)}{d\eta}\right)=(0,0)}$ and ${\displaystyle \left(u^*(c),\frac{d u^*(c)}{d\eta}\right)=(1,0)}$.  In particular, one can show \cite{KeenerSneyd} that
\[
u^*(c)(\eta)=\frac{1}{2}\left(1+\tanh\left(\frac{\eta}{2\sqrt{2}}\right)\right)
\]
is such a heteroclinic solution to (\ref{bisode}) with
\begin{equation}
c=\frac{1}{\sqrt{2}}(1-2\alpha).
\label{cval}
\end{equation}
Using the $\SEone$-equivariance, it follows that for any $\gamma\in\mathbb{R}$,
\begin{equation}
u(x,t)=u^*(c)(x+\gamma+ct)=\frac{1}{2}\left(1+\tanh\left(\frac{x+ct+\gamma}{2\sqrt{2}}\right)\right)
\label{twsol}
\end{equation}
is a travelling wave solution to (\ref{bis}) for $c$ given by (\ref{cval}).  These solutions are such that $u\rightarrow 1$ as $x\rightarrow\infty$ and $u\rightarrow 0$ as $x\rightarrow -\infty$, and represent a wave front which propagates from right to left when $c>0$ ($\alpha<1/2$) and a wave back propagating from left to right when $c <0$ ($\alpha>1/2$).
Note however that when $\alpha=1/2$, we have $c=0$ in (\ref{cval}), and thus (\ref{twsol}) represents a one-parameter family of steady-state solutions
\[
u(x,t)=u^*(0)(x+\gamma)=\frac{1}{2}\left(1+\tanh\left(\frac{x+\gamma}{2\sqrt{2}}\right)\right)\,\,\,\,\,\,\,\,\,\forall\gamma\in\mathbb{R}.
\]
The above analysis of (\ref{bis}) is standard for many classes of equations which admit travelling wave solutions in homogeneous media, i.e. suppose that the travelling-wave ansatz (\ref{twansatz}) holds, and seek homoclinic (for pulses) or heteroclinic (for fronts and backs) solutions in a reduced equation.  

Suppose now that we perturb (\ref{bis}) by introducing a small non-homogeneous term
\begin{equation}
\frac{\partial u}{\partial t}=\frac{\partial^2 u}{\partial x^2}+u(1-u)(u-\alpha)+\varepsilon\,g(x),
\label{bisp}
\end{equation}
where $g(x)$ is some smooth (non-constant) bounded function, and $0<\varepsilon\ll1$ is a small parameter.  Equation (\ref{bisp}) no longer admits the translation symmetry of (\ref{bis}), so pure travelling wave solutions such as (\ref{twansatz}) are unlikely to exist.  However, by continuous dependence of the solution on parameters, we expect solutions of (\ref{bisp}) which are ``close'' (in some sense) to (\ref{twsol}) for $\varepsilon$ small enough. 
There are essentially two regimes in parameter space to consider:

\noindent
{\em Regime I}: If the computed value of $c$ from formula
(\ref{cval}) is non-zero and large relative to $\varepsilon$, then heuristically speaking we expect to have a solution of (\ref{bisp}) which is approximately a wave front whose ``propagation speed'' is spatially-dependent ($c(x)\approx c$), but {\em non-vanishing on $\mathbb{R}$}.  In other words, in this regime, we expect travelling wave solutions of (\ref{bis}) to persist as solutions of (\ref{bisp}) with wave speed modulated by the heterogeneity, but the perturbation is not large enough to stop the wave from propagating along the real line. 

\noindent
{\em Regime II}: 
When $c$ and $\varepsilon$ are of the same order and are both small, propagation of the wave through the inhomogeneity may no longer be possible because of the possible persistence of steady states of (\ref{bis}) which may survive the perturbation in (\ref{bisp}).   This phenomenon is known as {\em propagation failure}, or {\em wave blocking}, and has been extensively studied in the literature for various specific models (e.g. \cite{BMR,IM,KFB,LewisKeener,TYBN,Xin}).

Our point of view is that most of the above discussion which was centered around the specific example (\ref{bis}) and (\ref{bisp}) does not depend qualitatively on the particular differential equation, but only on the fact that translation symmetry is broken by the addition of a small perturbation.  
We will thus consider a family of abstract differential equations parametrized by two parameters: $c$ and $\varepsilon$.  The parameter $\varepsilon$ controls the translational symmetry-breaking, while the parameter $c$ parametrizes the drift speed of relative equilibria (when $\varepsilon=0$).  Our main result will be to 
show that the point $(c,\varepsilon)=(0,0)$ in parameter space acts as an organizing center for the phenomenon of wave-blocking.  Indeed, we will show that wave-blocking occurs in a cone emanating from this organizing center,  and that this is a universal feature for any abstract differential equation which undergoes forced translational symmetry-breaking (under certain hypotheses).

To our knowledge, the study of the phenomenon of wave-blocking using techniques from dynamical systems for reaction-diffusion equations was first undertaken in \cite{LewisKeener}.  In that paper, the approach was to use a combination of comparison methods, bifurcation theory and numerical analysis to study a fixed travelling wave ($c>0$) in scalar reaction-diffusion equations, and let the size of the inhomogeneity increase until a saddle-node bifurcation occurs, and the resulting equilibria block propagation.  In this sense, the analysis in \cite{LewisKeener} is focussed away from the above-mentioned organizing center at $(c,\varepsilon)=(0,0)$.

Recently, in \cite{TYBN}, an analysis for a specific inhomogeneous system of reaction-diffusion equations was undertaken, using a formalism which is somewhat related to the approach we present here, i.e. centering the analysis around the codimension 2 drift bifurcation.  However, our approach will be general enough to include reaction-diffusion equations, as well as integro-differential equations, delay-reaction-diffusion equations, and potentially more.  
In \cite{KFB}, an analysis of wave-blocking in a specific system of integro-differential equations for inhomogeneous neural networks was undertaken, and equations for wave-blocking were derived which are similar in spirit to the equations we will derive here for general systems.

Our approach will be to use the center-manifold theorem of Sandstede, Scheel and Wulff \cite{SSW} for abstract group-equivariant evolution equations on Banach spaces as a key ingredient in our analysis.  Specifically, we consider the semi-linear equation
\begin{equation}
\frac{du}{dt}={\cal A}\,u+{\cal F}(u,c)+\varepsilon\,{\cal G}(u,\varepsilon)
\label{preode}
\end{equation}
which has $\SEone$-symmetry if and only if $\varepsilon=0$ (i.e. ${\cal G}$ breaks the translation symmetry).    When $\varepsilon=0$ this equation will admit a smooth two-dimensional invariant manifold of solutions foliated by travelling waves (relative equilibria) with drift speed $c\in (-c_0,c_0)$ for some $c_0>0$.   Assuming normal hyperbolicity, this manifold persists when $\varepsilon\neq 0$, but the dynamics on this manifold are perturbed.  In particular, the line of equilibria $c=0$ for the unperturbed equation is deformed, and it is this perturbed family of equilibria which causes propagation failure when $(c,\varepsilon)$ lies in some cone in parameter space which emanates from the point $(0,0)$.  Our main result to this effect is Theorem \ref{mainthm}.

The paper is organized as follows.  In section 2, we give the functional analytic group-equivariant framework for our analysis.  
In section 3 we consider the unperturbed (i.e. possessing full translation symmetry) differential equations, and apply the center manifold theorem of Sandstede, Scheel and Wulff \cite{SSW}.  The curve corresponding to zero drift speed on the center manifold will be central to the analysis of the perturbed equations, which is done in section 4 using the global parametrization near relative equilibria presented in \cite{SSW}.   In section 5, we briefly discuss modification of our method to study parameter-independent equations which admit families of travelling wave solutions such as the Fisher-KPP equation.  We end with some concluding remarks in section 6.

\Section{Functional analytic setup and translation symmetry}
Let $X$ be a Banach space with norm $||\,\,\,\,||$ and
consider the following nonlinear differential equation on $X$:
\begin{equation}
\frac{du}{dt}={\cal A}\,u+{\cal F}(u,c)+\varepsilon\,{\cal G}(u,\varepsilon)
\label{ode}
\end{equation}
where $c$ and $\varepsilon$ are real parameters, $\varepsilon$ small, ${\cal A}:X\longrightarrow X$ is a sectorial densely defined closed linear operator,  
the functions ${\cal F}:X^{\alpha}\times\mathbb{R}\longrightarrow X$ and ${\cal G}:X^{\alpha}\times\mathbb{R}\longrightarrow X$ are $C^{k+2}$ smooth for some $k\geq 1$ and some $\alpha\in [0,1)$ (see \cite{Henry}), 
and ${\cal G}$ is uniformly bounded on $X\times [0,\varepsilon_0]$ for some $\varepsilon_0>0$. 
Under these conditions, (\ref{ode}) generates a local semi-flow $\Phi_{t}(u;\varepsilon,c)$ on $Y:=X^{\alpha}$ \cite{Henry}.

Now, let $\Gamma=\SEone$ denote the 
group of rigid translations of the real line, which is isomorphic to the
additive group of real numbers.  Let
${\cal T}:\Gamma\longrightarrow GL(Y)$ be a faithful strongly continuous and isometric representation of $\Gamma$ in the space of bounded invertible operators on $Y$.
For $a\in\Gamma$, we will use the notation
${\cal T}_a$ to denote ${\cal T}(a)$.

The infinitesimal generator of ${\cal T}$ is the densely defined linear operator $\xi$, defined by 
\[
\xi\,u=\lim_{a\rightarrow 0}\,\frac{1}{a}({\cal T}_a-I)u,
\]
wherever this limit exists.

\begin{hyp}
We make the following assumptions on the functions ${\cal F}$ and ${\cal G}$ and the operators $\xi$ and ${\cal A}$:
\begin{enumerate}
\item[(a)] The domains of the linear operators ${\cal A}$ and $\xi$, ${\cal D}({\cal A})=Y_1$ and ${\cal D}(\xi)=Y_2$, are both dense and $\Gamma$-invariant.
Furthermore, there is a dense $\Gamma$-invariant subspace $Y_3\subset Y_1\cap Y_2$ such that ${\cal A}(Y_3)\subset Y_2$, $\xi(Y_3)\subset Y_1$, and
$\xi\,{\cal A}={\cal A}\,\xi$ on $Y_3$.  
\item[(b)] ${\cal A}{\cal T}_a={\cal T}_a{\cal A}$ on ${Y_1}$.
\item[(c)] The function ${\cal F}$ is $\Gamma$-equivariant, i.e. 
\begin{equation}
{\cal F}({\cal T}_au,c)={\cal T}_a{\cal F}(u,c),\,\,\,\,\,\forall u\in Y,\,c\in\mathbb{R}
\label{Fequiv}
\end{equation}
This implies also
\begin{equation}
D_u{\cal F}({\cal T}_au,c){\cal T}_av={\cal T}_aD_u{\cal F}(u,c)v,\,\,\,\,\,\forall\,u,v,\in Y,\,\,c\in\mathbb{R},
\label{DFequiv1}
\end{equation}
\begin{equation}
\xi\,{\cal F}(u,c)=D_u{\cal F}(u,c)\xi\,u\,\,\,\forall\,u\in Y_2,\,\, c\in\mathbb{R},\,\,\mbox{\rm such that}\,\,{\cal F}(u,c)\in Y_2,
\label{DFequiv3}
\end{equation}
and
\begin{equation}
D_c{\cal F}({\cal T}_au,c)={\cal T}_aD_c{\cal F}(u,c),\,\,\,\,\,\forall\,u\in Y,\,\,c\in\mathbb{R}
\label{DFequiv2}
\end{equation}
where $D_u{\cal F}$ and $D_c{\cal F}$ denote respectively the partial Fr\'echet derivatives of ${\cal F}$ with respect to $u$ and to $c$.
\item[(e)] The function ${\cal G}$ is not $\Gamma$-equivariant
\end{enumerate}
\label{hyp2}
\end{hyp}

We say that (\ref{ode}) undergoes {\em  forced translational symmetry-breaking} when $\varepsilon \neq 0$, since the equation admits $\Gamma$ as a symmetry group if and only if $\varepsilon=0$.

Here are a few examples of applications where Hypotheses \ref{hyp2} are satisfied.

%

\begin{examp}
Consider the following nonlinear reaction-diffusion equation
\begin{equation}
\frac{\partial u}{\partial t}(x,t)=D\frac{\partial^2u}{\partial x^2}(x,t)+f(u(x,t),c)+\varepsilon\,g(x),
\label{exam1}
\end{equation}
where $u(x)=(u_1,(x),\ldots,u_n(x))^T$, $D$ is a diagonal $n\times n$ matrix with positive entries,  $f:\mathbb{R}^n\times\mathbb{R}\longrightarrow\mathbb{R}^n$ is smooth, and $g:\mathbb{R}\longrightarrow\mathbb{R}^n$ is smooth, non-constant and bounded. 

We set $X={\cal L}^2(\mathbb{R},\mathbb{R}^n)$ and 
${\displaystyle Y_1={\cal D}\left(\frac{\partial^2}{\partial\,x^2}\right)=H^2(\mathbb{R},\mathbb{R}^{n})}$.
The action of $\Gamma=\SEone$ on $Y=X^\alpha$ is defined by
\[
({\cal T}_a\,u)(x)=u(x+a).
\]
The infinitesimal generator for this action is 
\[
(\xi\,u)(x)=u'(x)
\]
with dense $\Gamma$-invariant domain
\[
Y_2={\cal D}(\xi)=H^1(\mathbb{R},\mathbb{R}^{n}).
\]

Equation (\ref{exam1}) is invariant under this action if and only if $\varepsilon=0$.  
For the subspace $Y_3$ of Hypothesis \ref{hyp2} (a), we can choose
\[
Y_3=H^3(\mathbb{R},\mathbb{R}^{n})
\]
and can easily verify that ${\cal A}(Y_3)\subset Y_2$, $\xi(Y_3)\subset Y_1$ and ${\cal A}\,\xi=\xi\,{\cal A}$ on $Y_3$.
\label{examp1}
\end{examp}

\begin{examp}
Let ${\cal A}$ be a constant $n\times n$ matrix, and consider the following nonlinear integro-differential equation which arises frequently in neural field models with nonlocal interactions \cite{KFB}.
\begin{equation}
\frac{\partial u}{\partial t}(x,t)={\cal A}\,u(x,t)+\int_{-\infty}^\infty\,W(x-y)f(u(y,t),c)\,dy+\varepsilon\,\int_{-\infty}^{\infty}\,D(y)g(u(y,t))\,dy
\label{exam2}
\end{equation}
where
$W(\eta)$ and $D(\eta)$ are smooth $n\times n$ matrix-valued functions, $D$ is non-constant, and $f:\mathbb{R}^n\times\mathbb{R}\longrightarrow\mathbb{R}^n$ and 
$g:\mathbb{R}^n\longrightarrow\mathbb{R}^n$ are smooth and bounded.  Furthermore, $W$, $D$, $f$ and $g$ are such that the integrals in (\ref{exam2}) are all well-defined. 
We have $X=Y=C^0_{\mbox{\small unif}}(\mathbb{R},\mathbb{R}^n)$, and the $\SEone$-action is given by translation of the domain, as in Example \ref{examp1} above.   Equation (\ref{exam2}) has $\Gamma=\SEone$ symmetry if and only if $\varepsilon=0$.
We set
\[
Y_1=Y,\,\,\,\,\,Y_2=Y_3=C^1_{\mbox{\small unif}}(\mathbb{R},\mathbb{R}^{n}).
\]
\end{examp}

\begin{examp}
Let ${\cal H}$ be the Hilbert space ${\cal L}^2(\mathbb{R},\mathbb{R})$ and consider the Banach space $X=C^0([-r,0],{\cal H})$, $r>0$, with sup norm.  For a function $u:\mathbb{R}\longrightarrow {\cal H}$, define $u_t\in X$ by 
$u_t(\theta)(x)=u(t+\theta)(x)$ for $t\in\mathbb{R}$ and $\theta\in [-r,0]$.  Suppose $\Lambda:X\longrightarrow {\cal H}$ is a bounded linear operator, 
$A: {\cal D}(A)\subset {\cal H}\longrightarrow {\cal H}$ is the densely defined operator ${\displaystyle Au(t)=\frac{\partial^2v}{\partial x^2}(t,x)}$, where $u(t)(x)=v(t,x)$.
Consider the 
retarded partial functional differential equation
\begin{equation}
\dot{u}(t)=Au(t)+\Lambda(u_t)+F(u_t,c)+\varepsilon\,G(u_t)
\label{fde}
\end{equation}
where the nonlinearities $F:X\times\mathbb{R}\longrightarrow {\cal H}$ and $G : X\longrightarrow {\cal H}$ are sufficiently smooth.

The group $\SEone$ acts on $X$ using the standard action on ${\cal H}$ as in Example \ref{examp1} above:
\[
{\cal T}_a(u(\theta)(x))=u(\theta)(x+a).
\]
We suppose that $\Lambda$ and $F$ are such that
\[
\Lambda({\cal T}_a\,\varphi)(x)=\Lambda(\varphi)(x+a),\,\,\,\forall\,\varphi\in X,\,\forall\,a\in\mathbb{R}
\]
and
\[
F({\cal T}_a\,\varphi,c)(x)=F(\varphi,c)(x+a),\,\,\,\forall\,\varphi\in X,\,\forall\,a\in\mathbb{R},\,\forall\,c\in\mathbb{R},
\]
but $G$ does not satisfy this symmetry property.  We also assume that if $\varphi\in C^0([-r,0],H^1(\mathbb{R},\mathbb{R}))$, then $\Lambda(\varphi)\in H^1(\mathbb{R},\mathbb{R})$.  We note that these conditions on $\Lambda$ and $F$ are satisfied in many applications.

As was done in \cite{Hale}, we define 
\[
X_0=X_0(\theta)=\left\{\begin{array}{ccl}I_{\cal H}&\mbox{\rm if}&\theta=0\\0&\mbox{\rm if}&-r\leq\theta<0\end{array}\right.
\]
and consider the Banach space $BX$ of functions from $[-r,0]$ to $X$ which are uniformly continuous on $[-r,0)$ and with a jump discontinuity at $0$.  Elements of $BX$
can be written as $U=u+X_0\alpha$ with $u\in X$ and $\alpha\in {\cal H}$, so that $BX$ is identified with $X\times {\cal H}$.  The differential equation (\ref{fde}) can then be written as an abstract ODE on $BX$:
\begin{equation}
\frac{dv}{dt}={\cal A}v+X_0(F(v,c)+\varepsilon\,G(v,x))
\label{absfde}
\end{equation}
where ${\cal A}:BX\longrightarrow BX$ is defined by
\[
{\cal A}\,\psi=\dot{\psi}+X_0[\Lambda(\psi)+A\,\psi(0)-\dot{\psi}(0)],
\]
on the domain
\[
{\cal D}({\cal A})=
Y_1=\{\varphi\in X\,:\,\dot{\varphi}\in X,\,\varphi(0)\in H^2(\mathbb{R},\mathbb{R})\}.
\]

The action of $\SEone$ on $BX$ is:
\[
{\cal T}_a\,(u(\theta)(x)+X_0\alpha(x))=u(\theta)(x+a)+X_0\alpha(x+a).
\]
Equation (\ref{absfde}) has $\SEone$ symmetry if and only if $\varepsilon=0$.
The infinitesimal generator $\xi$ of $\SEone$ is 
\[
\xi\,(u(\theta)(x)+X_0\alpha(x))=u(\theta)'(x)+X_0\alpha'(x)
\]
defined on 
\[
{\cal D}(\xi)=Y_2=C^0([-r,0],H^1(\mathbb{R},\mathbb{R}))\times H^1(\mathbb{R},\mathbb{R})\subset X\times {\cal H}\cong BX,
\]
where the $'$ denotes a derivative with respect to $x$ for an element of $H^1(\mathbb{R},\mathbb{R})$.
We may then choose
\[
Y_3=\{\varphi\in X\,:\,\varphi,\,\dot{\varphi}\in C^0([-r,0],H^1(\mathbb{R},\mathbb{R})),\,\varphi(0)\in\,H^3(\mathbb{R},\mathbb{R})\}.
\]
It is then easy to verify that $Y_1$, $Y_2$ and $Y_3$ are $\SEone$-invariant, and that $\xi(Y_3)\subset Y_1$, ${\cal A}(Y_3)\subset Y_2$ and that $\xi\,{\cal A}={\cal A}\,\xi$ on $Y_3$.
\label{examp3}
\end{examp}

\Section{Relative equilibria in unperturbed equations}

For our analysis, it will be useful to single out the parameter $c$ in (\ref{ode}) and consider a suspended system
\begin{equation}
\begin{array}{lll}
{\displaystyle\frac{du}{dt}}&=&{\displaystyle{\cal A}\,u+{\cal F}(u,c)+\varepsilon\,{\cal G}(u,\varepsilon)}\\
&&\\
{\displaystyle\frac{dc}{dt}}&=&{\displaystyle 0}
\end{array}
\label{susode}
\end{equation}
on the Banach space $Y\times\mathbb{R}$, with norm $||\,\,\,\,\,||_{Y\times\mathbb{R}}$ defined by
\[
\left|\left|\left(\begin{array}{c}u\\c\end{array}\right)\right|\right|_{Y\times\mathbb{R}}=||u||+|c|.
\]
We extend the action ${\cal T}_a$ of $\Gamma$ to $Y\times \mathbb{R}$ using the trivial representation on the $\mathbb{R}$ factor, i.e.
\[
{\cal T}_a\left(\begin{array}{c}u\\c\end{array}\right)=\left(\begin{array}{c}{\cal T}_a\,u\\c\end{array}\right).
\]

If $\varepsilon=0$ in (\ref{susode}), then this equation is $\Gamma$-symmetric.  The $\Gamma$ symmetry is broken when $\varepsilon\neq 0$.
In this section, we will study (\ref{susode}) with $\varepsilon=0$, i.e.
\begin{equation}
\begin{array}{lll}
{\displaystyle\frac{du}{dt}}&=&{\displaystyle{\cal A}\,u+{\cal F}(u,c)}\\
&&\\
{\displaystyle\frac{dc}{dt}}&=&{\displaystyle 0}.
\end{array}
\label{presusode}
\end{equation}

For a given fixed value of $c$, a {\em relative equilibrium} for (\ref{presusode}) is a solution 
of the form
\[
\left(\begin{array}{c}u(t)\\c(t)\end{array}\right)=\left(\begin{array}{c}{\cal T}_{\gamma t}u^*\\c\end{array}\right),
\]
where $\gamma\in\mathbb{R}$, and
$u^*\in Y_3$ is such that ${\cal F}(u^*,c)\in Y_2$.
Substitution into the first equation of (\ref{susode}) and using equivariance yields
\begin{equation}
{\cal T}_{\gamma\,t}\gamma\xi\,u^*={\cal T}_{\gamma\,t}({\cal A}u^*+{\cal F}(u^*,c)).
\label{rel}
\end{equation}
If we assume that $\xi\,u^*\neq 0$, it follows that
\begin{equation}
|\gamma|=\frac{||{\cal A}u^*+{\cal F}(u^*,c)||}{||\xi\,u^*||}.
\label{propspeed}
\end{equation}
The real number $\gamma$ is called the {\em drift speed} of the relative equilibrium.
We now make the following hypothesis for (\ref{presusode}):
\begin{hyp}
We assume that there is an interval $(c_a,c_b)\subset\mathbb{R}$ and a smooth function 
\[
u^*:(c_a,c_b)\longrightarrow Y_3
\]
such that
\begin{enumerate}
\item[(a)] For all $c\in (c_a,c_b)$, we have ${\cal F}(u^*(c),c)\in Y_2$, $\xi\,u^*(c)\neq 0$, and both $\xi\,u^*(c)$ and ${\cal A}u^*(c)$ are smooth functions of $c$, with derivatives given respectively by $\xi\,u^*_c(c)$ and ${\cal A}u^*_c(c)$.  Furthermore, we assume that ${\cal T}_au^*(c)=u^*(c)$ if and only if $a=0$ (i.e. $u^*(c)$ has trivial isotropy).
\item[(b)] For all $c\in (c_a,c_b)$, (\ref{susode}) has a relative equilibrium 
\[
\left(\begin{array}{c}u(t)\\c(t)\end{array}\right)=\left(\begin{array}{c}{\cal T}_{\gamma(c)\,t}u^*(c)\\c\end{array}\right),
\]
 where the drift speed $\gamma(c)$ is the smooth function defined by
\begin{equation}
\gamma(c)\xi\,u^*(c)={\cal A}u^*(c)+{\cal F}(u^*(c),c).
\label{propspeed2}
\end{equation}
\item[(c)] The function $\gamma(c)$ defined above is monotone on $(c_a,c_b)$, and there exists $c_0\in (c_a,c_b)$ such that $\gamma(c_0)=0$.
\end{enumerate}
\label{uchyp}
\end{hyp}
\begin{prop}
Without loss of generality (upon changes of variables and relabelling of functions) we may assume that $\gamma(c)=c$ in (\ref{propspeed2}), i.e.
\begin{equation}
c\,\xi\,u^*(c)={\cal A}u^*(c)+{\cal F}(u^*(c),c),\,\,\,\,\mbox{\rm $\forall\,c$ near 0}.
\label{main_eigeneq}
\end{equation}
\end{prop}
\proof The function $\gamma(c)$ is invertible near $c=c_0$, so we may define $\tilde{c}=\gamma(c)$, $\tilde{u}^*(\tilde{c})=u^*(\gamma^{-1}(\tilde{c}))$ and
$\widetilde{\cal F}(u,\tilde{c})={\cal F}(u,\gamma^{-1}(\tilde{c}))$.  Equation (\ref{propspeed2}) then becomes
\[
\tilde{c}\,\xi\,u^*(\gamma^{-1}(\tilde{c}))={\cal A}u^*(\gamma^{-1}(\tilde{c}))+{\cal F}(u^*(\gamma^{-1}(\tilde{c})),\gamma^{-1}(\tilde{c})),
\]
or
\[
\tilde{c}\,\xi\,\tilde{u}^*(\tilde{c})={\cal A}\tilde{u}^*(\tilde{c})+\widetilde{\cal F}(\tilde{u}^*(\tilde{c}),\tilde{c}).
\]
We get (\ref{main_eigeneq}) by relabelling $\tilde{c}\rightarrow c$, $\tilde{u^*}\rightarrow u^*$ and $\widetilde{\cal F}\rightarrow {\cal F}$. 
\qed

 In the sequel, we will assume that this change of parameter and relabelling of functions has already been made in (\ref{presusode}) and (\ref{susode}).

Setting $c=0$ in (\ref{main_eigeneq}) and applying $\xi$ yields
\begin{equation}
{\cal A}\,\xi\,u^*(0)+D_u{\cal F}(u^*(0),0)\xi\,u^*(0)=0
\label{1stevect}
\end{equation}
Also, differentiation of (\ref{main_eigeneq}) with respect to $c$ and evaluating at $c=0$ yields
\begin{equation}
{\cal A}\,u^*_c(0)+D_u{\cal F}(u^*(0),0)\,u^*_c(0)+D_c{\cal F}(u^*(0),0)=\xi\,u^*(0)
\label{2ndevect}
\end{equation}
Thus, the following densely defined closed linear operator
\begin{equation}
\begin{array}{c}
{\displaystyle L:Y\times\mathbb{R}\longrightarrow Y\times\mathbb{R}}\\
\\
L\left(\begin{array}{c}\varphi\\\omega\end{array}\right)=\left(
\begin{array}{c}
{\cal A}\,\varphi+D_u{\cal F}(u^*(0),0)\,\varphi+D_c{\cal F}(u^*(0),0)\,\omega\\0\end{array}\right)
\end{array}
\label{Ldef}
\end{equation}
is such that
\[
L\left(\begin{array}{c}u^*_c(0)\\1\end{array}\right)=\left(\begin{array}{c}\xi\,u^*(0)\\0\end{array}\right)\,\,\,\,\,\,\,\,\,\mbox{\rm and}\,\,\,\,\,\,\,\,\,\,
L\left(\begin{array}{c}\xi\,u^*(0)\\0\end{array}\right)=\left(\begin{array}{c}0\\0\end{array}\right),
\]
i.e. $L$ has an eigenvalue at $0$ with eigenvector $\phi_1\equiv (\xi\,u^*(0),0)^T$ and generalized eigenvector $\phi_2\equiv (u^*_c(0),1)^T$.   We make the following further assumptions:
\begin{hyp}
We will assume that 
\begin{enumerate}
\item[(a)] The eigenvalue at $0$ of ${L}$ is such that  the corresponding generalized eigenspace $E$ is two-dimensional, i.e.
$E=\mbox{\rm span}\{\phi_1,\phi_2\}$.  Furthermore, all other elements of the spectrum of ${L}$ are bounded away from the imaginary axis (spectral gap condition).
We assume that the Fredholm index of $L$ is 0, and that we have the following $L$-invariant splitting
\begin{equation}
Y\times\mathbb{R}=E\oplus W=\mbox{\rm span}\{\phi_1,\phi_2\}\oplus W
\label{splitY1}
\end{equation}
with bounded projection $P$, $\mbox{\rm range}(P)=E$, $\ker\,(P)=W$, and $L(W)=W$.  We will also use the projection
\[
Q:E\longrightarrow E
\]
such that $\mbox{\rm range}(Q)=\mathbb{R}\{\phi_2\}$.   Note that since any element of $W=L(W)$ must have the form
\[
\left(\begin{array}{c}\mu\\0\end{array}\right),\,\,\,\mu\in Y
\]
and since ${\displaystyle \phi_1=\left(\begin{array}{c}\xi\,u^*(0)\\0\end{array}\right)}$,
then 
\[
QP\left(\begin{array}{c}\sigma\\0\end{array}\right)=\left(\begin{array}{c}0\\0\end{array}\right),\,\,\forall\,\sigma\in Y.
\]
\item[(b)] The map ${\displaystyle a\mapsto {\cal T}_au^*(0)}$ is $C^{k+2}$ in $a\in\Gamma$.
\item[(c)] For any $\kappa>0$, there exists a $\delta > 0 $ such that $||{\cal T}_au^*(0)-u^*(0)||\geq\delta$ whenever $|a|\geq\kappa$.
\item[(d)] The map $a\mapsto {\cal T}_av$ is $C^{k+1}$ whenever $v\in E$.
\item[(e)] If $P$ denotes the spectral projection onto $E$ in item (a) above, then the projections ${\cal T}_aP{\cal T}_{-a}$ are $C^{k+1}$ smooth in $a\in\Gamma$ in the operator norm.
\end{enumerate}
\label{spectralhyp}
\end{hyp}
\begin{prop}
The mapping $u^*(c)$ in Hypothesis \ref{uchyp} can be chosen so that
\begin{equation}
P\left(\begin{array}{c}u^*(c)-u^*(0)\\0\end{array}\right)=\left(\begin{array}{c}0\\0\end{array}\right),\,\,\,\,\mbox{\rm for all $c$ near $0$}.
\label{orth}
\end{equation}
\label{orthprop}
\end{prop}
\proof
Clearly ${\displaystyle QP\left(\begin{array}{c}u^*(c)-u^*(0)\\0\end{array}\right)=\left(\begin{array}{c}0\\0\end{array}\right)}$ (see the end of item (a) in Hypothesis \ref{spectralhyp}).  Consider the smooth map
\[
f:\mathbb{R}\times\mathbb{R}\longrightarrow \mathbb{R}\{\phi_1\}
\]
defined by
\[
f(a,c)=(I-Q)P\left(\begin{array}{c}{\cal T}_{-a}\,u^*(c)-u^*(0)\\0\end{array}\right).
\]
We have $f(0,0)=0$ and
\[
f_a(0,0)=-(I-Q)P\phi_1=-\phi_1\neq 0.
\]
By the implicit function theorem, there exists a smooth function $a(c)$ such that $a(0)=0$ and such that
\[
(I-Q)P\left(\begin{array}{c}{\cal T}_{-a(c)}\,u^*(c)-u^*(0)\\0\end{array}\right)=\left(\begin{array}{c}0\\0\end{array}\right),\,\,\,\,\,\mbox{\rm for all $c$ near $0$}.
\]
We then redefine ${\cal T}_{-a(c)}\,u^*(c)\longrightarrow u^*(c)$.\qed
\begin{rmk}
It follows from Proposition \ref{orthprop} that 
\[
P\left(\begin{array}{c}u^*_c(0)\\0\end{array}\right)=\left(\begin{array}{c}0\\0\end{array}\right).
\]
\label{rmk1}
\end{rmk}

Hypothesis \ref{spectralhyp} allows us to use the center manifold theorem of Sandstede, Scheel and Wulff \cite{SSW}, from which it follows that the set
\[
{\cal S}=\left\{\left(\begin{array}{c}{\cal T}_au^*(c)\\c\end{array}\right)\,\,;\,\,a\in\Gamma, \mbox{\rm $c$ near 0}\right\}
\]
is a normally hyperbolic $C^{k+1}$ smooth invariant center manifold for (\ref{presusode}).  The dynamics of (\ref{presusode}) restricted to ${\cal S}$ is
rather trivial, described by the two-dimensional center manifold ordinary differential equations
\begin{equation}
\begin{array}{rcl}
\dot{a}&=&c\\
\dot{c}&=&0.
\end{array}
\label{cmeqs_0}
\end{equation}
In particular, the line of equilibria $c=0$ for (\ref{cmeqs_0}) corresponds to the smooth one-dimensional group orbit 
\[
\Theta=\left\{\,\left(\begin{array}{c}{\cal T}_a\,u^*(0)\\0\end{array}\right)\,;\,a\in\Gamma\,\right\}
\]
of the steady state ${\displaystyle\left(\begin{array}{c}u^*(0)\\0\end{array}\right)}$ of (\ref{presusode}).  

\begin{figure}[htbp]
\begin{center}
\input{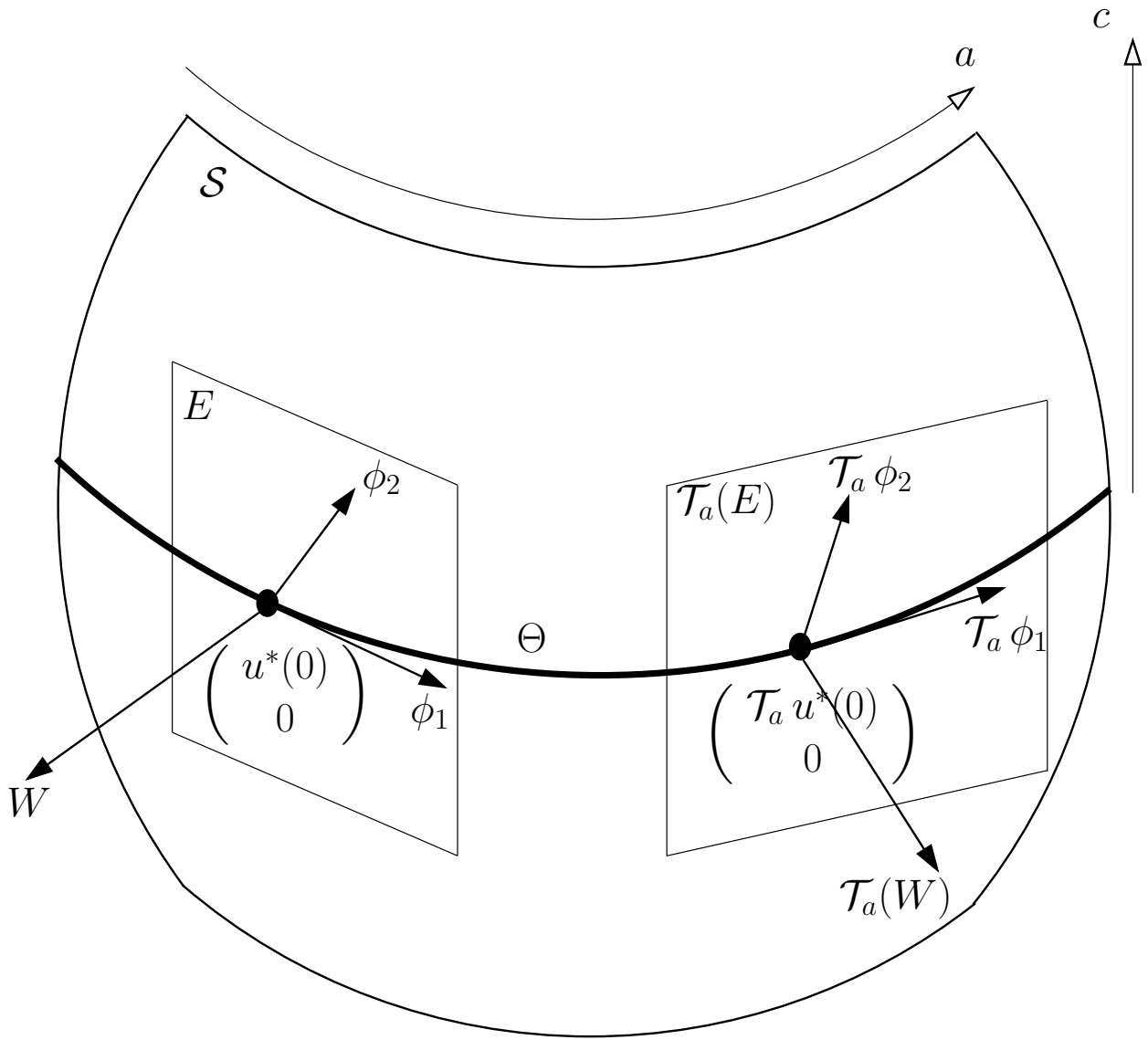}
\caption{{Schematic representation of the decomposition of Hypothesis \ref{spectralhyp}(a) as it relates to the center manifold ${\cal S}$, the family of steady-states $\Theta$, and the coordinates $(a,c)$ of the center manifold equations (\ref{cmeqs_0})}.}
\label{fig1}
\end{center}
\end{figure}

\begin{figure}[htbp]
\begin{center}
\input{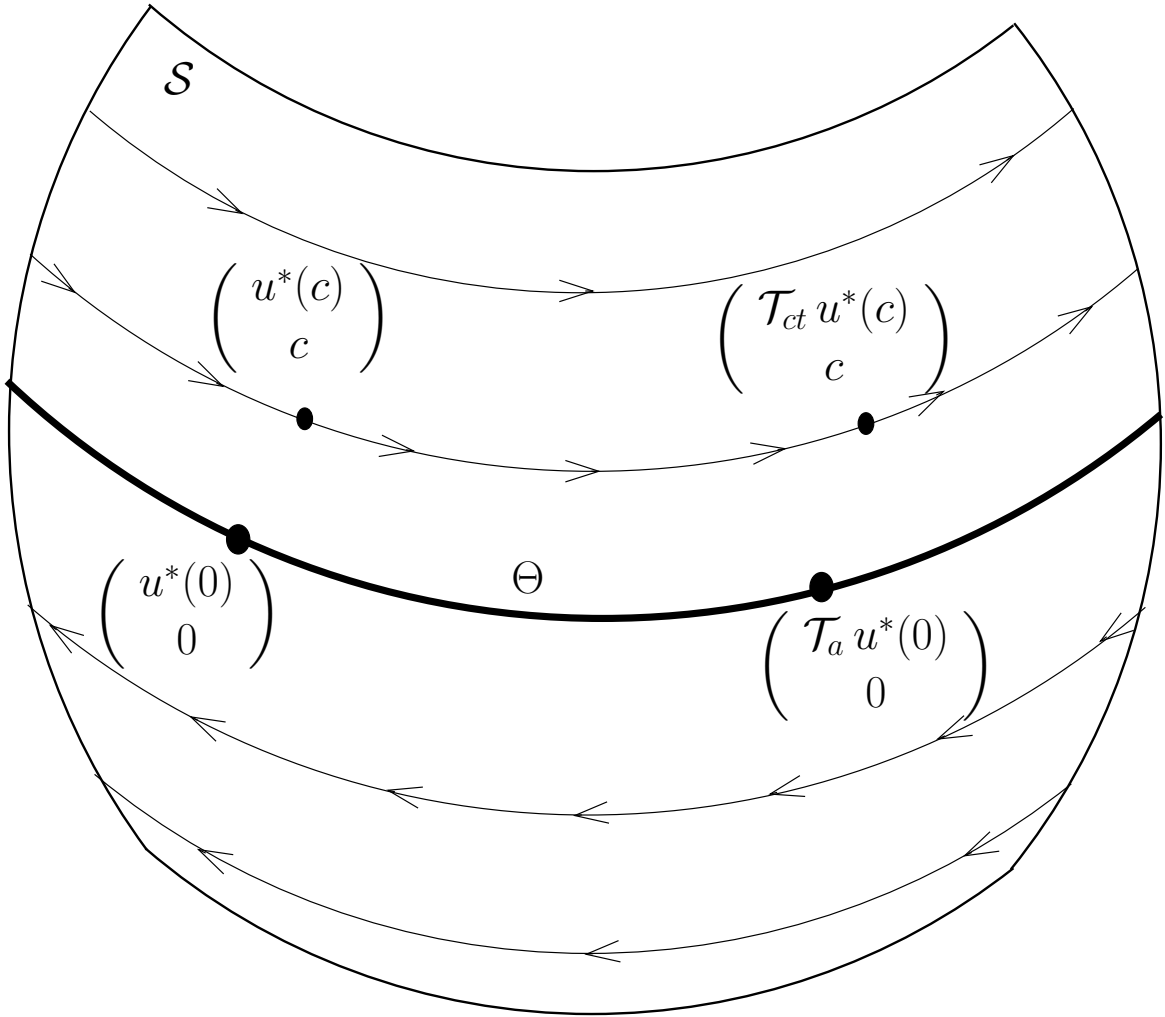}
\caption{Projection of the flow of (\ref{cmeqs_0}) onto the center manifold ${\cal S}$ for $\varepsilon=0$.  Closer arrows indicate slower flow.  The curve $\Theta$ ($c=0$) is a curve of equilibria ($\dot{a}=0$).}
\label{fig2}
\end{center}
\end{figure}

As was shown in \cite{SSW}, we obtain a global parametrization of a neighborhood of $\Theta$ as follows: any point $u$ close enough to $\Theta$ can be written as
\[
u=U_1(u)+U_2(u)+w(u),
\]
where
\[
U_1(u)={\cal T}_{a(u)}\left(\begin{array}{c}u^*(0)\\0\end{array}\right),\,\,\,\,U_2(u)\in {\cal T}_{a(u)}(\mbox{\rm range}(QP)),\,\,\,\,
w(u)\in {\cal T}_{a(u)}(W),
\]
where all functions of $u$ above are $C^{k+1}$.  In particular, for any $c$ close enough to $0$, we may write the following representation of the center manifold
${\cal S}$:
\begin{equation}
\left(\begin{array}{c}{\cal T}_au^*(c)\\c\end{array}\right)=
{\cal T}_a\left[\left(\begin{array}{c}u^*(0)\\0\end{array}\right)+c\phi_2+\left(\begin{array}{c}v(c)\\0\end{array}\right)\right],
\label{split}
\end{equation}
where by Proposition \ref{orthprop} and Remark \ref{rmk1} we have
\[
{\displaystyle
\left(\begin{array}{c}v(c)\\0\end{array}\right)=\left(\begin{array}{c}u^*(c)-u^*(0)-cu^*_c(0)\\0\end{array}\right) \in W
}, 
\]
i.e.
\[
P\left(\begin{array}{c}v(c)\\0\end{array}\right)=\left(\begin{array}{c}0\\0\end{array}\right),
\]
where $P$ is the projection operator onto $E$ in Hypothesis \ref{spectralhyp}(a).  See figures \ref{fig1} and \ref{fig2} for a schematic representation of 
${\cal S}$ and of the flow of (\ref{cmeqs_0}) on ${\cal S}$.

Taking into account the equations (\ref{cmeqs_0}) on the center manifold, the decomposition (\ref{split}) when substituted into (\ref{susode}) leads to the compatibility condition
\begin{equation}
c^2\,\,P\left(\begin{array}{c}\xi\,u^*_c(0)\\0\end{array}\right)+c\,P\left(\begin{array}{c}\xi\,v(c)\\0\end{array}\right)=
P\left(\begin{array}{c}H(c)\\0\end{array}\right)
\label{compat}
\end{equation}
where
\[
\begin{array}{l}
H(c)={\cal F}(u^*(0)+cu^*_c(0)+v(c),c)-{\cal F}(u^*(0),0)-c\,D_u{\cal F}(u^*(0),0)u^*_c(0)\\
\\
-c\,D_c{\cal F}(u^*(0),0)-D_u{\cal F}(u^*(0),0)v(c).
\end{array}
\]

Using normal hyperbolicity arguments similar to those of LeBlanc and Wulff \cite{LW}, the set ${\cal S}$ will persist as a normally hyperbolic invariant manifold ${\cal S}_{\varepsilon}\cong {\cal S}$ for (\ref{susode}) for all $\varepsilon$ sufficiently close to $0$, but the dynamics on ${\cal S}_{\varepsilon}$ may be different than on ${\cal S}$.  In particular, we are interested in the persistence of equilibrium points on ${\cal S}_{\varepsilon}$ near $\Theta$.  

\Section{Translational symmetry-breaking}

We now consider equations (\ref{susode}) in the case where $\varepsilon\neq 0$.
Suppose ${\displaystyle\left(\begin{array}{c}u(t;c,\varepsilon)\\c\end{array}\right)}$ is a solution of (\ref{susode}) on the manifold ${\cal S}_{\varepsilon}$ near $\Theta$, for $0<\varepsilon \ll 1$.
We use the local coordinates (\ref{split}) and write
\begin{equation}
\left(\begin{array}{c}u(t;c,\varepsilon)\\c\end{array}\right)=
{\cal T}_{a(t)}\left[\left(\begin{array}{c}u^*(0)\\0\end{array}\right)+c\phi_2+\left(\begin{array}{c}v(c)\\0\end{array}\right)+\varepsilon\,\left(\begin{array}{c}E_1(a(t),c,\varepsilon)\\0\end{array}\right)\right]
\label{paraS}
\end{equation}
where $E=(E_1,0)^T$ is a smooth function 
\[
E:\mathbb{R}\times [-c_0,c_0]\times [0,\varepsilon_0]\longrightarrow (Y_3\times\mathbb{R})\cap W
\]
satisfying
\[
\sup_{a\in\mathbb{R}}\,\{||E(a,c,\varepsilon)||_{Y\times\mathbb{R}}\}\leq K<\infty\,\,\,\forall\, (c,\varepsilon)\,\in\,[-c_0,c_0]\times [0,\varepsilon_0]
\]
for some $c_0>0$, $\varepsilon_0>0$.
Substitution of (\ref{paraS}) into (\ref{susode}) yields
\[
\begin{array}{l}
\dot{a}\,\left(\phi_1+c\,\left(\begin{array}{c}\xi\,u^*_c(0)\\0\end{array}\right)+\left(\begin{array}{c}\xi\,v(c)\\0\end{array}\right)+O(\varepsilon)\right)=
c\,\phi_1+L\left(\begin{array}{c}v(c)\\0\end{array}\right)+\varepsilon\,L\left(\begin{array}{c}E_1\\0\end{array}\right)+\\
\\
\varepsilon\,\left(\begin{array}{c}{\cal T}_{-a}{\cal G}({\cal T}_a\,u^*(0),0)\\0\end{array}\right)
+\left(\begin{array}{c}{\cal F}(u^*(0)+c\,u^*_c(0)+v(c),c)\\0\end{array}\right)
-\left(\begin{array}{c}{\cal F}(u^*(0),0)\\0\end{array}\right)\\
\\-c\,\left(\begin{array}{c}D_u{\cal F}(u^*(0),0)u^*_c(0)\\0\end{array}\right)
-\left(\begin{array}{c}D_u{\cal F}(u^*(0),0)\,v(c)\\0\end{array}\right)-
\left(\begin{array}{c}D_c{\cal F}(u^*(0),0)\,c\\0\end{array}\right)
+\varepsilon\,\left(\begin{array}{c}\tilde{q}(a,c,\varepsilon)\\0\end{array}\right),
\end{array}
\]
where 
\[
\begin{array}{l}
\tilde{q}(a,c,\varepsilon)=\\
\\
{\displaystyle\frac{1}{\varepsilon}[{\cal F}(u^*(0)+c\,u^*_c(0)+v+\varepsilon\,E_1,c)-
{\cal F}(u^*(0)+c\,u^*_c(0)+v,c)]}
-D_u{\cal F}(u^*(0)+c\,u^*_c(0)+v,c)E_1\\
\\
+{\displaystyle D_u{\cal F}(u^*(0)+c\,u^*_c(0)+v,c)E_1-D_u{\cal F}(u^*(0),0)E_1}+\\
\\
+{\cal T}_{-a}{\cal G}({\cal T}_{a}(u^*(0)+c\,u^*_c(0)+v_1+\varepsilon\,E_1),\varepsilon)-{\cal T}_{-a}{\cal G}({\cal T}_a\,u^*(0),0)
\end{array}
\] 
is smooth and such that ${\displaystyle \lim_{(c,\varepsilon)\rightarrow (0,0)}\tilde{q}(a,c,\varepsilon)=0}$.   
Projecting the above equation onto $E$ using the projection $P$ and using (\ref{compat}) leads to
the following perturbation of the center manifold equations (\ref{cmeqs_0})
\begin{equation}
\begin{array}{rcl}
\dot{a}&=&c-\varepsilon\,r(a)+\varepsilon\,q(a,c,\varepsilon)\\
\dot{c}&=&0
\end{array}
\label{maincmeq}
\end{equation}
where
\begin{equation}
r(a)\phi_1=
-(I-Q)P\left(\begin{array}{c}{\cal T}_{-a}{\cal G}({\cal T}_a\,u^*(0),0)\\0\end{array}\right)
\label{rlambdadef}
\end{equation}
and $q(a,c,\varepsilon)$ is smooth and such that ${\displaystyle \lim_{(c,\varepsilon)\rightarrow (0,0)}\,q(a,c,\varepsilon)=0}$.

\begin{rmk}
We make the following remarks concerning (\ref{maincmeq}) and (\ref{rlambdadef}):
\begin{enumerate}
\item[(1)] If ${\cal G}$ were $\Gamma$ equivariant, then it would follow that $r(a)$ is a constant.  Since it is assumed that ${\cal G}$ breaks translational symmetry, then in general $r$ is a smooth, non-constant function of $a$.
\item[(2)] The function
$r(a)$ and the term $q(a,c,\varepsilon)$ in (\ref{maincmeq}) are uniformly bounded in $a$ for $(c,\varepsilon)$ near $(0,0)$.
\end{enumerate}
\end{rmk}
\begin{prop}
For any $\varepsilon$ sufficiently close to 0, there exists a smooth curve $\Theta_{\varepsilon}$ which is the graph of a function $c=\varepsilon\,r(a)+\varepsilon\,\sigma(a,\varepsilon)$ (where $\sigma$ is bounded and $\sigma(a,0)=0)$ such that $\dot{a}=0$ in (\ref{maincmeq}) whenever $(a,c)$ belongs to this curve.  Consequently, $\Theta_{\varepsilon}$ designates a curve of equilibria for the perturbed center manifold equations (\ref{maincmeq}).
\end{prop}
\proof This is a simple application of the implicit function theorem.
\qed

From this proposition, we conclude that if $c$ satisfies the relation
\[
\inf_{a\in\mathbb{R}} \{r(a)+\sigma(a,\varepsilon)\}<\frac{c}{\varepsilon}<\sup_{a\in\mathbb{R}} \{r(a)+\sigma(a,\varepsilon)\}
\]
then $\dot{a}$ in (\ref{maincmeq}) will become 0 at at least one point $a\in\mathbb{R}$.  Consequently, the flow line $(a(t),c)$ for (\ref{maincmeq}) contains at least one equilibrium point, which prevents the flow line from going from $a=-\infty$ to $a=\infty$.  The physical interpretation of this phenomenon is what is commonly referred to in the literature as {\em propagation failure}, or {\em wave-blocking} \cite{BMR,IM,KFB,LewisKeener,TYBN,Xin}.

\begin{figure}[htbp]
\begin{center}
\input{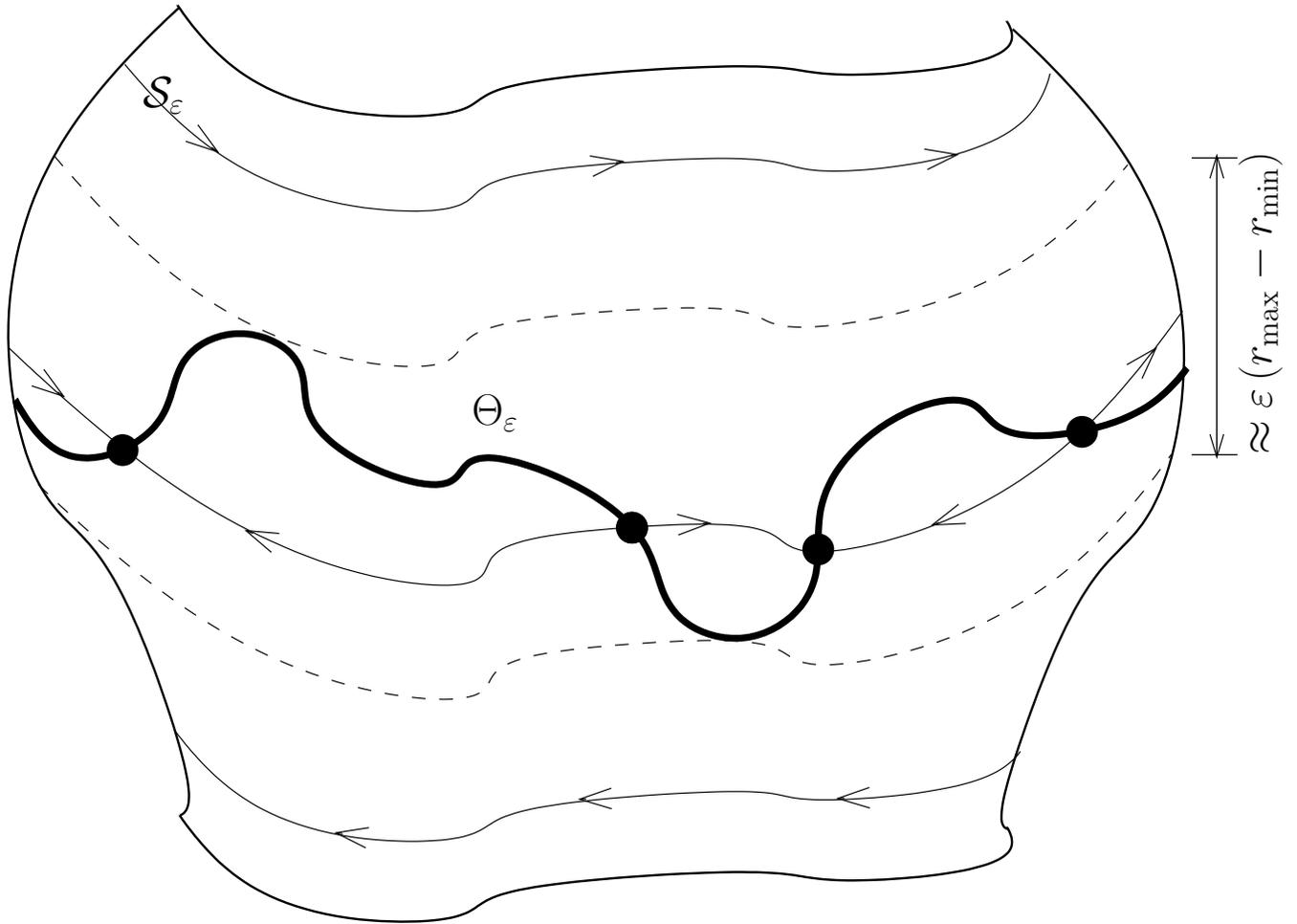}
\caption{Flow of (\ref{maincmeq}) on the perturbed center manifold ${\cal S}_{\varepsilon}$.  All flow lines have constant value of $c$, but $\dot{a}$ is not constant along the flow lines.  Flow lines between the dashed lines undergo wave-blocking (by meeting equilibrium points belonging to $\Theta_{\varepsilon}$).  If $r_{\max}$ and $r_{\min}$ designate respectively the supremum and the infimum of the function $r(a)$ on $\mathbb{R}$, then the width of $c$-values for which wave-blocking occurs is approximately $\varepsilon\,(r_{\max}-r_{\min})$.}
\label{fig3}
\end{center}
\end{figure}

On the other hand, if $c$ satisfies the relation
\[
\frac{c}{\varepsilon}>\sup_{a\in\mathbb{R}} \{r(a)+\sigma(a,\varepsilon)\}\,\,\,\,\,\,\mbox{\rm or}\,\,\,\,\,\,\,
\frac{c}{\varepsilon}<\inf_{a\in\mathbb{R}} \{r(a)+\sigma(a,\varepsilon)\}
\]
then $\dot{a}$ in (\ref{maincmeq}) will be strictly positive, or strictly negative, but non-constant.  Therefore the corresponding flow line $(a(t),c)$ for (\ref{maincmeq}) will evolve at a non-constant drift speed, but the sign of this drift speed will not change on the flow line.  Therefore, propagation failure (or wave-blocking) does not occur for this flow line.  The situation is depicted in figure \ref{fig3}.

\begin{figure}[htbp]
\begin{center}
\input{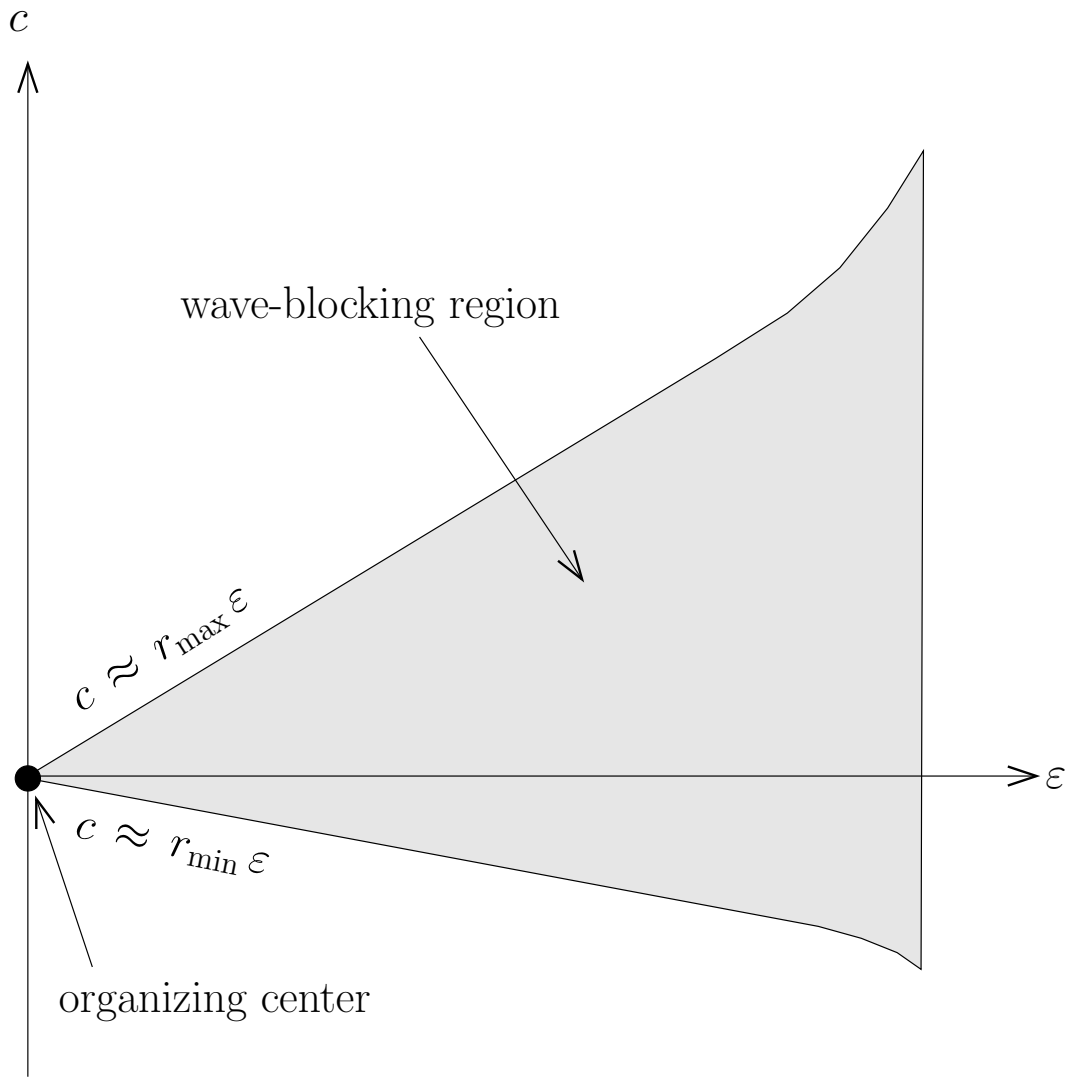}
\caption{Wave blocking occurs in (\ref{susode}) for parameter values $(c,\varepsilon)$ belonging to the illustrated cone.}
\label{fig4}
\end{center}
\end{figure}

We have thus proved the following
\begin{thm}
Consider the abstract differential equation (\ref{susode}) satisfying Hypotheses \ref{hyp2}, \ref{uchyp} and \ref{spectralhyp}.  For all $\varepsilon\geq 0$ close enough to $0$, the semi-flow for this differential equation admits a two-dimensional normally hyperbolic
$C^{k+1}$-smooth invariant manifold ${\cal S}_{\varepsilon}$.  
Solutions of (\ref{susode}) which belong to ${\cal S}_{0}\equiv {\cal S}$ correspond to relative equilibria (with constant drift speed $c$).  On the other hand, 
solutions of (\ref{susode}) which belong to ${\cal S}_{\varepsilon}$ (for $\varepsilon\neq 0$) correspond to perturbed relative equilibria, with non-constant drift speed.  In the latter case, propagation failure (or wave-blocking) occurs if $(c,\varepsilon)$ belongs to a cone in parameter space, as illustrated in figure \ref{fig4}.
\label{mainthm}
\end{thm}

\Section{Parameter independent cases}

Hypothesis \ref{spectralhyp}(a) implies that $D_c{\cal F}(u^*(0),0)\neq 0$, otherwise $L$ would also have $(0,1)^T$ as an eigenvector.
Thus, our results do not immediately apply to equations for which $D_c{\cal F}(u^*(0),0)= 0$, which includes the cases where ${\cal F}$ does not
depend explicitly on the parameter $c$, or on any other parameters.  This comment is particularly relevant for the Fisher-KPP equation \cite{Murray}:
\[
\frac{\partial u}{\partial t}=\frac{\partial^2u}{\partial x^2}+u(1-u)
\]
which admits travelling wave front solutions $U(x+ct)$ for all $c$ in a neighbourhood of $0$\footnote{For $-2<c<2$, the wave front has an oscillatory approach to 0.  Therefore, for many applications in biology, solutions with $-2<c<2$ are not considered since they are not biologically relevant (the waveform attains negative values as it oscillates).}, even though the reaction term $u(1-u)$ does not depend on any parameters.  For cases such as this, the main ideas in this paper can still be used, but need to be slightly modified.  

First, the equation 
\begin{equation}
\frac{du}{dt}={\cal A}\,u+{\cal F}(u)+\varepsilon\,{\cal G}(u,\varepsilon)
\label{nop}
\end{equation}
is studied without suspending the system (as was done in (\ref{susode})).   We then assume hypotheses similar to Hypotheses \ref{uchyp}, i.e. existence of a family of travelling wave solutions for (\ref{nop}) (when $\varepsilon=0$) parametrized by drift speed $c$.

The linear operator $L$ becomes
\[
\begin{array}{c}
L:Y\longrightarrow Y\\
\\
L(\varphi)={\cal A}\varphi+D_u{\cal F}(u^*(0))\varphi
\end{array}
\]
with eigenvector $\phi_1=\xi\,u^*(0)$ and generalized eigenvector $\phi_2=u^*_c(0)$.  The rest of the analysis proceeds essentially the same as in section 3 and 4.
In particular, one can show that there exists  a normally hyperbolic two-dimensional center manifold ${\cal S}=\{\,{\cal T}_a\,u^*(c)\,;\,a\in\mathbb{R},c\in (-c_0,c_0)\,\}\subset Y$ on which the flow of (\ref{nop}) at $\varepsilon=0$ is given by (\ref{cmeqs_0}), and which persists to a normally hyperbolic invariant manifold ${\cal S}_{\varepsilon}\cong {\cal S}$ for (\ref{nop}) when $\varepsilon\neq 0$.  However, the perturbed equations on ${\cal S}_{\varepsilon}$ become
\begin{equation}
\begin{array}{rcl}
\dot{a}&=&c-\varepsilon\,r_1(a)+\varepsilon\,q_1(a,c,\varepsilon)\\
\dot{c}&=&\varepsilon\,r_2(a)+\varepsilon\,q_2(a,c,\varepsilon)
\end{array}
\label{maincmeq2}
\end{equation}
where
\begin{equation}
r_1(a)\phi_1=
-(I-Q)P\left(\begin{array}{c}{\cal T}_{-a}{\cal G}({\cal T}_a\,u^*(0),0)\\0\end{array}\right),\,\,\,\,\,
r_2(a)\phi_2=QP\left(\begin{array}{c}{\cal T}_{-a}{\cal G}({\cal T}_a\,u^*(0),0)\\0\end{array}\right),
\label{rlambdadef2}
\end{equation}
and $q_{1,2}(a,c,\varepsilon)$ are smooth and such that ${\displaystyle \lim_{(c,\varepsilon)\rightarrow (0,0)}\,q_{1,2}(a,c,\varepsilon)=0}$.  

\begin{figure}[htbp]
\begin{center}
\input{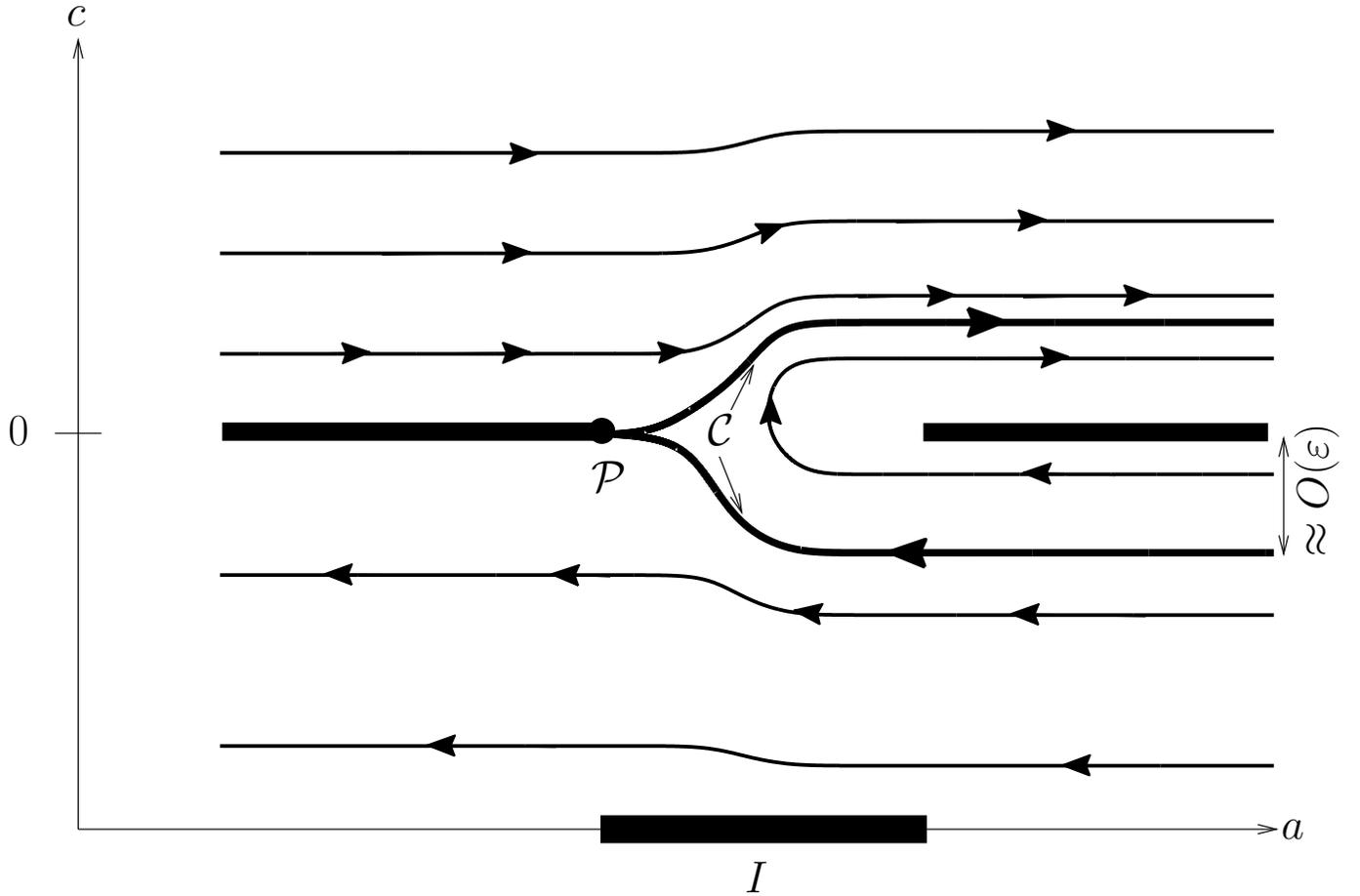}
\caption{Flow of (\ref{maincmeq3}) in the case described in Example \ref{inh21}.  Orbits on one side of ${\cal C}$ do not undergo propagation failure, however their drift speed is modulated by the inhomogeneity.  On the other side of ${\cal C}$, orbits coming from $a=\infty$ in a region of width $\approx O(\varepsilon)$ undergo propagation failure: they get reflected by the inhomogeneity and bounce back to $a=\infty$.  The two wide lines at $c=0$ on either sides of the interval $I$ correspond to unperturbed steady-states $c=0$ in (\ref{cmeqs_0}).}
\label{fig5}
\end{center}
\end{figure}

\begin{figure}[htbp]
\begin{center}
\input{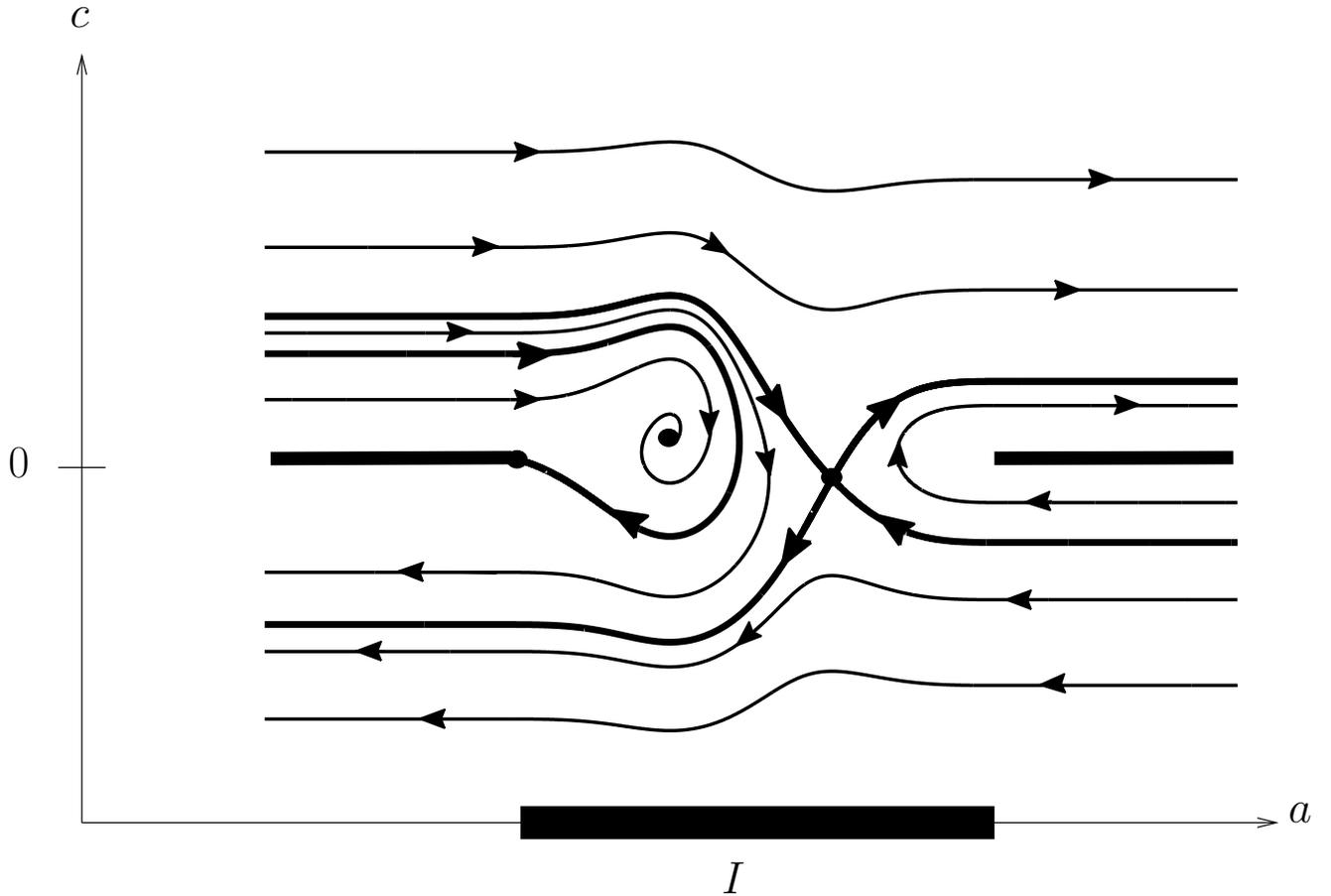}
\caption{Flow of (\ref{maincmeq3}) in the case described in Example \ref{inh22}.  
The stable and unstable manifolds of the saddle separate the phase space into 3 regions.  In two of these regions, orbits do not undergo propagation failure, however their drift speed is modulated by the inhomogeneity.  In the other region (of width $O(\varepsilon)$), orbits either reflect off the inhomogeneity, or get attracted into the stable node.}
\label{fig6}
\end{center}
\end{figure}

\begin{figure}[htbp]
\begin{center}
\input{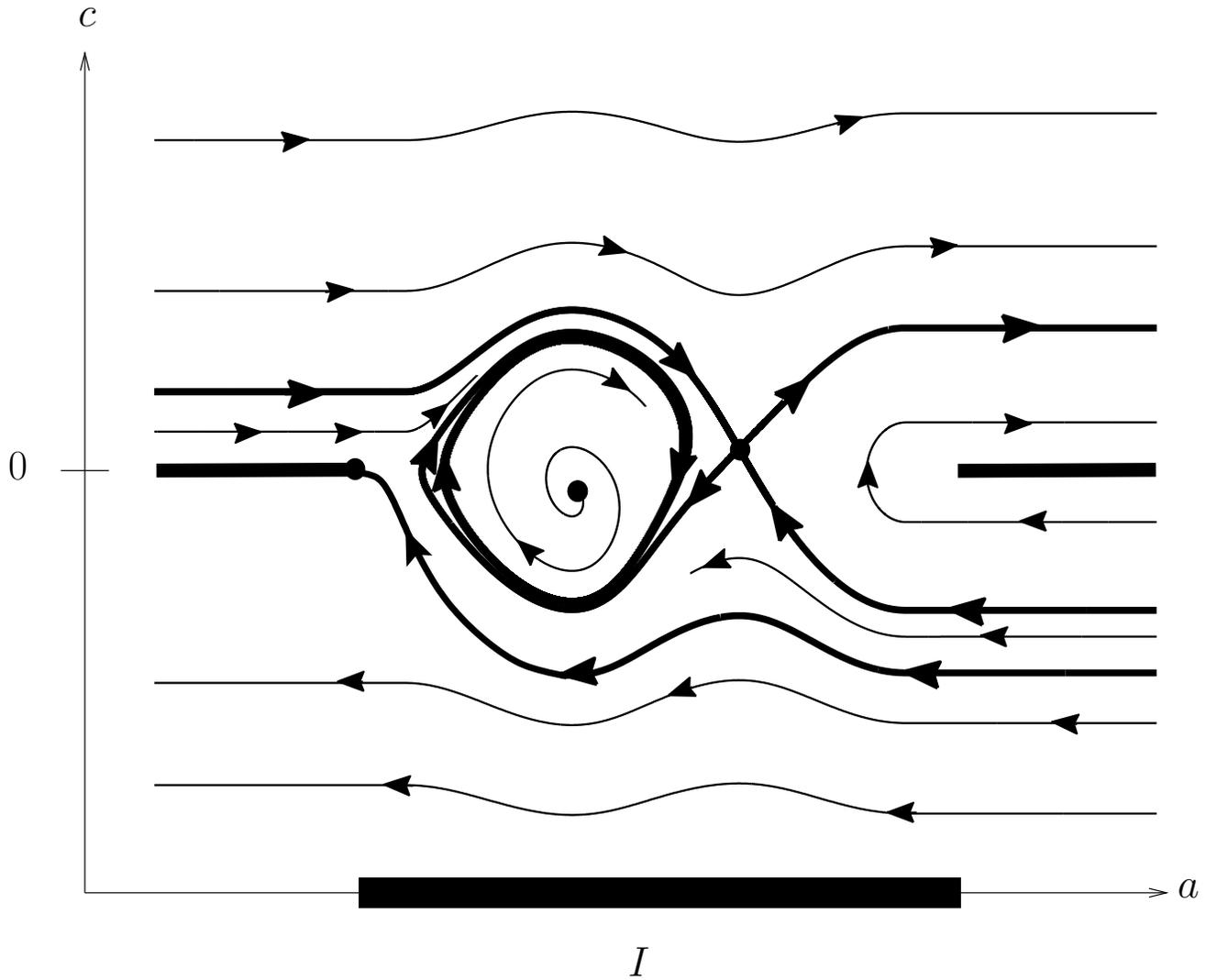}
\caption{Flow of (\ref{maincmeq3}) in the case described in Example \ref{inh23}.  Similar to the previous case, the stable and unstable manifolds of the saddle separate the phase space into regions in which orbits propagate along the entire real line, and regions where propagation failure occurs.  In this case, a stable limit cycle attracts certain orbits coming from $a=-\infty$.}
\label{fig7}
\end{center}
\end{figure}

For $(c,\varepsilon)$ near $(0,0)$, the principal part of (\ref{maincmeq2}) is 
\begin{equation}
\begin{array}{rcl}
\dot{a}&=&c-\varepsilon\,r_1(a)\\
\dot{c}&=&\varepsilon\,r_2(a).
\end{array}
\label{maincmeq3}
\end{equation}
The dynamics of (\ref{maincmeq3}) can be quite complicated, depending on the functions $r_1(a)$ and $r_2(a)$.  We will not give an exhaustive classification,
but we will give a few illustrative examples.   Many more examples are given in \cite{TYBN} for a specific three-species system of reaction-diffusion equations.
In order to simplify things, we will suppose in our examples that the inhomogeneity is localized in the sense that both $r_1(a)$ and $r_2(a)$ vanish outside the interval $I=(a_0,a_1)$.
\begin{examp}
Suppose the function $r_2(a)$ is of the same sign on $I$ (without loss of generality, suppose $r_2(a)>0$).  Then the phase space is separated into two regions by a curve ${\cal C}$ which is the union of two orbits of (\ref{maincmeq3}): one which tends to ${\cal P}=(a,c)=(a_0,\varepsilon\,r_1(a_0))$ as $t\rightarrow\infty$, and the other orbit which tends to this same point as $t\rightarrow-\infty$, as illustrated in figure \ref{fig5}.
On one side of ${\cal C}$, orbits with $c\neq 0$ propagate from $a=-\infty$ to $a=\infty$ (or vice-versa) without blocking, however the drift speed of the corresponding travelling wave is modulated by the inhomogeneity.  On the other side of ${\cal C}$, orbits of (\ref{maincmeq3}) undergo a reflection by the inhomogeneity as they travel from $a=\infty$ to the inhomogeneity, and then back to $a=\infty$, as illustrated in figure \ref{fig5}.  The size of the interval of $c$-values for which orbits are reflected is approximately $O(\varepsilon)$.
\label{inh21}
\end{examp}
\begin{examp}
Suppose now that the function $r_2(a)$ changes sign over $I$.  For the purposes of this example, we will suppose that $r_2(a)$ has two simple zeroes at $a=\rho_1$ and $a=\rho_2$ in $(a_0,a_1)$, and suppose that $r_2(a)>0$ on $(a_0,\rho_1)$ and on $(\rho_2,a_2)$, and $r_2(a)<0$ on $(\rho_1,\rho_2)$.  Then (\ref{maincmeq3}) has hyperbolic equilibrium points at $(\rho_1,\varepsilon\,r_1(\rho_1))$ (node) and at $(\rho_2,\varepsilon\,r_1(\rho_2))$ (saddle).  A possible phase diagram is shown in figure \ref{fig6}.  The stable and unstable manifolds of the saddle point separate the phase space into three regions.   In two of these regions, orbits propagate from $a=-\infty$ to $a=\infty$ (or vice-versa) without blocking, whereas in the other region, propagation failure occurs by one of two mechanisms: either reflection off the inhomogeneity, or blocking by the stable node.
\label{inh22}
\end{examp}
\begin{examp}
In Example \ref{inh22} above, instead if having a stable node, the node could be unstable and surrounded by a stable limit cycle.   The situation is depicted in figure \ref{fig7}.  This illustrates another mechanism for propagation failure: the travelling wave tends asymptotically to a state characterized by a waveform which oscillates about a steady state.
\label{inh23}
\end{examp}

As mentioned above, these three example only serve to illustrate some of the possibilities for the perturbed dynamics in (\ref{maincmeq3}).  An exhaustive classification would be quite tedious, and certainly beyond the scope of this paper.

\Section{Conclusion and Discussion}

We have presented a model-independent analysis of the phenomenon of wave-blocking, showing that under some general hypotheses, this phenomenon is universal in the context of forced translational symmetry-breaking in parametrized families of differential equations on Banach spaces.  The approach is general enough in scope to be applicable to reaction-diffusion partial differential equations, integro-differential equations, retarded partial functional differential equations, and more.
Our main theorem to that effect is Theorem \ref{mainthm}.

We note that from a practical point of view, explicit computations of the various terms in the perturbed center manifold equations (\ref{maincmeq}) and (\ref{maincmeq2}) require explicit knowledge of the various terms which contribute to it, for example the function $u^*(c)$, and the eigenvectors $\phi_1$ and $\phi_2$ of the linearization $L$.

We have chosen to work in a Banach space setting so that our results could apply to as many situations as possible.  In certain applications (e.g. when studying pulses which decay to 0 in reaction-diffusion partial differential equations), it is customary to work in the Hilbert space ${\cal L}^2(\mathbb{R},\mathbb{R}^n)$ (e.g. \cite{KFB, TYBN}).  In this case, one can use the adjoint operator $L^*$, the Fredholm alternative, and orthogonal projections for the various projections which were introduced in Section 3.

Finally, we have discussed in Section 5 how our approach can be modified to study parameter-independent equations (such as the Fisher-KPP equation) which admit families of travelling waves parametrized by drift speed.  We gave a few illustrative examples of how localized inhomogeneities can effect the dynamics of travelling waves in these cases.

\vspace*{0.25in}
\noindent
{\Large\bf Acknowledgments}

\vspace*{0.2in}
This research is partly supported by the
Natural Sciences and Engineering Research Council of Canada in the
form of a Discovery Grant (VGL) and a Postgraduate Scholarship (CR).



\begin{thebibliography}{10}

\bibitem{AA}
N.~Akhmediev and A.~Ankiewicz.
\newblock {\em Nonlinear pulses and beams.}
\newblock Chapman and Hall, London, (1997).

\bibitem{BNPR}
H.~Berestycki, G.~Nadin, B.~Perthame and L.~Ryzhik.
\newblock The non-local Fisher-KPP equation: travelling waves and steady states.
\newblock {\em Nonlinearity} {\bf 22}, (2009) 2813--2844.

\bibitem{BMR}
Y.~Boubendir, V.~M\'{e}ndez and H.G.~Rotstein.
\newblock Dynamics of one- and two-dimensional fronts in a bistable equation with time-delayed global
feedback: Propagation failure and control mechanisms.
\newblock {\em Phys. Rev. E.} {\bf 82}, (2010) 036601-1--20.

\bibitem{Ei}
S.-I.~Ei.
\newblock The Motion of Weakly Interacting Pulses in Reaction-Diffusion Systems.
\newblock {\em J. Dyn. Diff. Eqs.} {\bf 14}, (2002) 85--137.

\bibitem{FT}
T.~Faria and S.~Trofimchuk.
\newblock Positive travelling fronts for reaction-diffusion systems with distributed delay.
\newblock {\em Nonlinearity} {\bf 23}, (2010), 2457--2481.

\bibitem{FM}
P.C.~Fife and J.B.~Mcleod.
\newblock The approach of solutions of nonlinear diffusion
equations to travelling front solutions. 
\newblock {\em Arch. Rat. Mech. Anal.} {\bf 65}, (1977), 335--361.

\bibitem{GXZY}
Q.~Gan, R.~Xu, X.~Zhang and P.~Yang.
\newblock Travelling waves of a three-species Lotka-Volterra food-chain model with spatial diffusion and time delays.
\newblock {\em Nonlinear Anal. Real World Appl.} {\bf 11}, (2010), 2817--2832.

\bibitem{GSS2}
M.~Golubitsky, I.~Stewart and D.G.~Schaeffer.
\newblock {\em Singularities and Groups in Bifurcation Theory, Vol. II.}
\newblock Applied Mathematical Sciences {\bf 69}, Springer-Verlag, New York, (1988).

\bibitem{GK}
G.A.~Gottwald and L.~Kramer.
\newblock On propagation failure in one- and two-dimensional excitable media.
\newblock {\em Chaos} {\bf 14}, (2004) 855--863.

\bibitem{Hale}
J.K.~Hale.
\newblock Critical cases for neutral functional differential equations.
\newblock {\em J. Diff. Eqs.} {\bf 10}, (1971) 59--82.

\bibitem{Henry}
D.~Henry. 
\newblock {\em Geometric theory of semilinear parabolic equations.}
\newblock Lecture Notes in Mathematics {\bf 804}, Springer-Verlag, New York, (1981).

\bibitem{IM}
H.~Ikeda and M.~Mimura.
\newblock Wave-blocking phenomena in bistable reaction-diffusion systems. 
\newblock {\em SIAM. J. Appl. Math.} {\bf 49}, (1989) 515-538.

\bibitem{Jones}
C.K.R.T.~Jones.
\newblock Stability of the travelling wave solution of the FitzHugh-Nagumo system. 
\newblock {\em Trans. AMS} {\bf  286}, (1984) 431--469.

\bibitem{KeenerSneyd}
J.~Keener and J.~Sneyd.
\newblock
{\em Mathematical Physiology.}
\newblock Interdisciplinary Applied Mathematics {\bf 8}, Springer-Verlag,
New York, (1998).

\bibitem{KFB}
Z.P.~Kilpatrick, S.E.~Folias and P.C.~Bressloff.
\newblock travelling Pulses and Wave Propagation Failure in Inhomogeneous Neural Media.
\newblock {\em SIAM J. App. Dyn. Syst.} {\bf 7}, (2008) 161--185.

\bibitem{LW}
V.G.~LeBlanc and C.~Wulff
\newblock Translational symmetry-breaking for spiral waves.
\newblock {\em J. Nonlin. Sci.} {\bf 10}, (2000) 569--601.

\bibitem{LewisKeener}
T.J.~Lewis and J.P. Keener.
\newblock Wave-block in excitable media due to regions of depressed excitability.
\newblock {\em SIAM. J. Appl. Math.} {\bf 61}, (2000) 293--316.

\bibitem{Murray}
J.D.~Murray.
\newblock
{\em Mathematical Biology: 1. An Introduction, 3rd ed.}
\newblock Interdisciplinary Applied Mathematics {\bf 17}, Springer-Verlag, New York, (2002).

\bibitem{ReyMackey}
A.D.~Rey and M.C.~Mackey.
\newblock Bifurcations and travelling waves in a delayed partial differential equation.
\newblock {\em Chaos} {\bf 2}, (1992) 231--244.

\bibitem{Sandstede}
B.~Sandstede.
\newblock Stability of travelling waves.
\newblock In: {\em Handbook of Dynamical Systems II} (B Fiedler, ed.).
North-Holland (2002) 983-1055.

\bibitem{SSW}
B.~Sandstede, A.~Scheel and C.~Wulff.
\newblock Dynamics of spiral waves on unbounded domains using center-manifold reductions.
\newblock {\em J. Diff. Eqs} {\bf 141} (1997) 122--149.

\bibitem{SRLC}
Z.~Szyma\'{n}ska, C.M.~Rodrigo, M.~Lachowicz and M.A.J.~Chaplain.
\newblock Mathematical modelling of cancer invasion of tissue: the role and effect of nonlocal interactions.
\newblock {\em Math. Models Methods Appl. Sci.} {\bf 19}, (2009) 257--281.

\bibitem{TYBN}
T.~Teramoto, X.~Yuan, M.~B\"{a}r and Y.~Nishiura.
\newblock Onset of unidirectional pulse propagation in an excitable medium with asymmetric heterogeneity.
\newblock {\em Phys. Rev. E.} {\bf 79} (2009) 046205-1--16.

\bibitem{VV}
A.I.~Volpert, V.A.~Volpert and V.A.~Volpert.
\newblock
{\em travelling waves solutions of parabolic systems.}
\newblock Transl. Math. Mono. {\bf 140}, Amer. Math. Soc., Providence, (1994).

\bibitem{Xin}
J.~Xin.
\newblock Front Propagation in Heterogeneous Media.
\newblock {\em SIAM Rev.} {\bf 42} (2000) 161--230.





\end{thebibliography}
\end{document}